\documentclass[12pt,fleqn]{article}
\usepackage{latexsym,amsfonts,amssymb,amsmath,amsthm,longtable,epsfig,graphics,lscape}
\usepackage[cp1251]{inputenc}
\usepackage[english]{babel}

\newcommand{\ba} {\begin{array}}
\newcommand{\ea} {\end{array}}
\newcommand{\be} {\begin{equation}}
\newcommand{\ee} {\end{equation}}

{\theoremstyle{definition}

}

\begin{document}

\begin{center}
{\bf Complex dynamic of the system of nonlinear difference equations in the Hilbert space.}
\end{center}

{\bf Pokutnyi O.O. Institute of mathematics of NAS of Ukraine} \\

{\it {\bf Abstract.} In the given article the necessary and sufficient conditions of the existence of solutions of boundary value problem for the nonlinear system  in the Hilbert spaces are obtained. Examples of such systems like a Lotka-Volterra are considered.  Bifurcation and branching conditions of solutions are obtained. }

{\bf Introduction.} The system of difference equations is the subject of numerous publications, and it is impossible to analyze all of them in detail. In this article we develop constructive methods of analysis of linear and weakly nonlinear boundary-value problems for difference equations, which occupy a central place in the qualitative theory of dynamical systems. We consider such problems that the operator of the linear part of the equation does not have an inverse. Such problems include the so called critical (or resonance) problems (when considering problem can have nonunique solution and not for any right hand sides). We use the well-known technique of generalized inverse operators \cite{BoiSam} and use the notion of a strong generalized solution of an operator equation developed in \cite{LyashNomir}. In this way, one can prove the existence of solutions of different types for the system of operator equations in the Hilbert spaces. There exist three possible types of solutions: classical solutions, strong generalized solutions, and strong pseudosolutions \cite{Tikhonov}. For the analysis of a weakly nonlinear system, we develop the well-known Lyapunov-Schmidt method. This approach gives possibility to investigate  many problems in difference equations and mathematical biology from a single point of view.

{\bf Statement of the problem.} Consider the following boundary value problem
\begin{equation} \label{eq:1}
x(n + 1, \varepsilon) = a(n)x(n, \varepsilon) + b(n)y(n, \varepsilon) + \varepsilon Z_{1}(x(n, \varepsilon), y(n, \varepsilon), n, \varepsilon) + f_{1}(n);
\end{equation}
\begin{equation}\label{eq:2}
y(n + 1, \varepsilon) = c(n)x(n, \varepsilon) + d(n)y(n, \varepsilon) + \varepsilon Z_{2}(x(n, \varepsilon), y(n, \varepsilon), n, \varepsilon) + f_{2}(n);
\end{equation}
\begin{equation} \label{eq:3}
l \left( \begin{array}{ccc} x(\cdot, \varepsilon) \\
 y(\cdot, \varepsilon) \end{array} \right) = \alpha,
\end{equation}
where operators  $ \{ a(n), b(n), c(n), d(n) \in \mathcal{L}(\mathcal{H}), n \in J \subset \mathbb{Z} \}$, $\mathcal{L}(\mathcal{H})$  is the space of linear and bounded operators which acts from $\mathcal{H}$ into itself,  vector-functions $f_{1}(n), f_{2}(n) \in l_{\infty}(J, \mathcal{H})$,
$$
l_{\infty}(J, \mathcal{H}) = \{ f : J \rightarrow \mathcal{H}, ||f||_{l_{\infty}} = \sup_{n \in J}||f(n)||_{\mathcal{H}} < \infty \},
$$
$Z_{1}, Z_{2}$ are smooth nonlinearities; a linear and bounded operator $l$ translates solutions of (\ref{eq:1}),
(\ref{eq:2}) into the Hilbert space $\mathcal{H}_{1}$, $\alpha$ is an element of the space $\mathcal{H}_{1}$, $\alpha \in \mathcal{H}_{1}$.
(instead of $l_{\infty}(J, \mathcal{H})$ we can consider another functional space $\mathcal{T}(J, \mathcal{L}(\mathcal{H}))$).
We find solutions of the boundary value problem (\ref{eq:1})-(\ref{eq:3}) which for $\varepsilon = 0$ turns in one of solutions of generating boundary value problem
\begin{equation} \label{eq:4}
x_{0}(n + 1) = a(n)x_{0}(n) + b(n)y_{0}(n) + f_{1}(n);
\end{equation}
\begin{equation} \label{eq:5}
y_{0}(n+1) = c(n)x_{0}(n) + d(n)y_{0}(n) + f_{2}(n);
\end{equation}
\begin{equation}\label{eq:6}
l\left( \begin{array}{ccc}x_{0}(\cdot) \\
y_{0}(\cdot) \end{array} \right) = \alpha.
\end{equation}

{\bf Linear case.} Consider the following vector $z_{0}(n) = (x_{0}(n), y_{0}(n))^{T}$, sequence of operator matrices
$$
A_{n} = \left(\begin{array}{ccc} a(n) & b(n) \\
c(n) & d(n)
\end{array} \right),
$$
and sequence of vector-functions $f(n) = (f_{1}(n), f_{2}(n))^{T}$. Here $T$ is a transpose operation. Then we can rewrite the generating boundary value problem (\ref{eq:4})-(\ref{eq:6}) in the following  form
\begin{equation} \label{eq:7}
z_{0}(n + 1) = A_{n}z_{0}(n) + f(n),
\end{equation}
\begin{equation} \label{eq:8}
lz_{0}(\cdot) = \alpha.
\end{equation}
Define an operator $\Phi(m, n) = A_{m + 1}A_{m}... A_{n + 1}, m > n $, $\Phi(m, m) = I$. The operator $U(m)  = \Phi(m, 0)$ is an evolution operator \cite{Chuesh}.
General solution $z_{0}(n)$ of (\ref{eq:7}) can be represented in the following form:
\begin{equation} \label{equat:100}
z_{0}(n) = \Phi(n, 0)z_{0} + g(n),
\end{equation}
where
$$
g(n) = \sum_{i  = 0}^{n}\Phi(n, i)f(i).
$$

{\bf Remark 1.} {\it It should be noted that if the sequence of operator matrices $A_{n}$ has bounded inverse $A_{n}^{-1} \in \mathcal{L}(\mathcal{H})$, then the general solution of (\ref{eq:7}) can be represented in the following form
$$
z_{0}(n) = U(n)z_{0} + \sum_{i = 0}^{n} U(n)U^{-1}(i)f(i).
$$
}
Substituting  representation (\ref{equat:100}) in the boundary condition (\ref{eq:8}) we obtain the following operator equation
\begin{equation} \label{equat:101}
Qz_{0} = h,
\end{equation}
where the operator $Q$ and the element $h$ have the following form
$$
Q = l\Phi(\cdot, 0),~~~~Q: \mathcal{H} \rightarrow \mathcal{H}_{1},  ~~~ h = \alpha - lg(\cdot).
$$
According to the theory of generalized solutions which was represented in \cite{BoiPok} and theory of Moore-Penrose pseudoinvertible operators \cite{BoiSam} for the equation (\ref{equat:101}) we have the following variants:

1) Suppose that $R(Q) = \overline{R(Q)}$ ($R(Q)$ is the image of the operator $Q$). In this case we have that the equation (\ref{equat:101}) is solvable if and only if the following
condition is hold \cite{BoiSam}:
\begin{equation} \label{equat:105}
\mathcal{P}_{N(Q^{*})}h = 0.
\end{equation}
Here $\mathcal{P}_{N(Q^{*})}$ is an orthoprojector onto the kernel of adjoint operator $Q^{*}$ ($\mathcal{P}_{N(Q^{*})} = \mathcal{P}_{N(Q^{*})}^{2} = \mathcal{P}_{N(Q^{*})}^{*}$).
Under condition (\ref{equat:105}) the set of solutions of (\ref{equat:101}) has the following form \cite{BoiSam}:
$$
z_{0} = Q^{+}h + \mathcal{P}_{N(Q)}c, ~~\forall c \in \mathcal{H},
$$
where $Q^{+}$ is Moore-Penrose pseudoinverse   \cite{BoiSam}, \cite{Moor}, \cite{Penr} to the operator $Q$, $\mathcal{P}_{N(Q)}$ is orthoprojector onto the kernel of the operator $Q$.

2) Consider the case when $R(Q) \neq \overline{R(Q)}$. In this case there is strong Moore-Penrose pseudoinverse $\overline{Q}^{+}$  \cite{BoiPok} to the operator $\overline{Q}$ ($\overline{Q} : \overline{\mathcal{H}} \rightarrow \mathcal{H}_{1}$ is extension of the operator $Q$ onto extended space $\mathcal{H} \subset \overline{\mathcal{H}}$ \cite{BoiPok}). Condition of generalized solvability has the following form:
\begin{equation} \label{equat:106}
\mathcal{P}_{N(\overline{Q}^{*})}h = 0.
\end{equation}
Condition (\ref{equat:106}) guarantees only that $h \in \overline{R(Q)}$.
Under condition (\ref{equat:106}) the set of generalized solutions of the equation (\ref{equat:101}) has the following form:
\begin{equation}\label{equat:107}
z_{0} = \overline{Q}^{+}h + \mathcal{P}_{N(\overline{Q})}c, ~~\forall c \in \mathcal{H}.
\end{equation}
If $h \in R(Q)$ then generalized solutions will be classical.

3) Suppose that $R(Q) \neq \overline{R(Q)}$ and $h \notin \overline{R(Q)}$. It means that the following condition is hold
\begin{equation} \label{equat:106}
\mathcal{P}_{N(\overline{Q}^{*})}h \neq 0.
\end{equation}
Under condition (\ref{equat:106}) the set of generalized quasisolutions \cite{BoiPok}, \cite{BoiSam} has the following form:
$$
z_{0} = \overline{Q}^{+}h + \mathcal{P}_{N(\overline{Q})}c, ~~\forall c \in \mathcal{H}.
$$
In such a way we obtain the following theorem.

{\bf Theorem 1.} {\it Boundary value problem (\ref{eq:7}), (\ref{eq:8}) is solvable.

a1) There are generalized solutions of (\ref{eq:7}), (\ref{eq:8}) if and only if
\begin{equation} \label{eq:9}
\mathcal{P}_{N(\overline{Q}^{*})} \{ \alpha - l\sum_{i = 0}^{\cdot}\Phi(\cdot, i)f(i) \} = 0,
\end{equation}
where the operator $\mathcal{P}_{N(\overline{Q}^{*})}$ is an orthoprojector onto the kernel of operator $\overline{Q}^{*}$ ($\overline{Q}$ is the extended operator to the operator $Q$). If the element $(\alpha - l\sum_{i = 0}^{\cdot}\Phi(\cdot, i)f(i)) \in R(Q)$ then solutions will be classical solutions;

b1)   under condition (\ref{eq:9}) the set of generalized solutions of the boundary value problem (\ref{eq:7}), (\ref{eq:8}) has the following form
$$
z_{0}(n, c) = \overline{G[f, \alpha]}(n) + \mathcal{P}_{N(\overline{Q})}c, \forall c \in \mathcal{H},
$$
where the generalized Green operator has the form
$$
\overline{G[f, \alpha]}(n)  = \Phi(n, 0)\overline{Q}^{+}\{ \alpha - l\sum_{i = 0}^{\cdot} \Phi(\cdot, i)f(i)\};
$$
a2) There are strong quasisolutions of (\ref{eq:7}), (\ref{eq:8}) if and only if the following condition is hold
\begin{equation} \label{eq:10}
\mathcal{P}_{N(\overline{Q}^{*})}\{ \alpha - l\sum_{i = 0}^{\cdot}\Phi(\cdot, i)f(i) \} \neq 0;
\end{equation}

b1)   under condition (\ref{eq:10}) the set of strong quasisolutions of the boundary value problem (\ref{eq:7}), (\ref{eq:8}) has the following form
$$
z_{0}(n, c) = \overline{G[f,\alpha]}(n) +\mathcal{P}_{N(\overline{Q})}c, \forall c \in \mathcal{H}.
$$
}

{\bf Nonlinear case.}
 Consider the nonlinear boundary value problem (\ref{eq:1})-(\ref{eq:3}).  Use denotions we can rewrite this problem in the following form
\begin{equation} \label{eq:11}
z(n + 1, \varepsilon) = A_{n}z(n, \varepsilon) + \varepsilon Z(z^{T}(n, \varepsilon), n, \varepsilon),
\end{equation}
\begin{equation} \label{eq:12}
lz(\cdot, \varepsilon) = \alpha.
\end{equation}

{\bf Theorem 2.} (necessary condition). {\it Suppose that the boundary value problem (\ref{eq:11}), (\ref{eq:12}) has solution
$z(n, \varepsilon)$ which for $\varepsilon = 0$ turns in one of solutions $z_{0}(n, c)$ with element $c \in \mathcal{H}$ ($z(n, 0) = z_{0}(n, c)$). Then $c$ satisfies the following operator equation for generating elements
\begin{equation} \label{eq:13}
F(c) = \mathcal{P}_{N(\overline{Q}^{*})}l\sum_{i = 0}^{\cdot}\Phi(\cdot, i)Z(z_{0}^{T}(i, c), i, 0) =
\end{equation}
$$
 = \mathcal{P}_{N(\overline{Q}^{*})}l\sum_{i = 0}^{\cdot}\Phi(\cdot, i)Z(\overline{G[f, \alpha]}(i) + \mathcal{P}_{N(\overline{Q})}c, \cdot, 0) = 0.
$$
}
{\bf Proof.} According to the theorem 1  the boundary value problem (\ref{eq:11}), (\ref{eq:12}) has solution if and only if the following condition is true:
\begin{equation} \label{eq:14}
\mathcal{P}_{N(\overline{Q}^{*})}\{\alpha - l\sum_{i = 0}^{\cdot}\Phi(\cdot, i)(f(i) + \varepsilon Z(z^{T}(i, \varepsilon), i, \varepsilon)) \} = 0.
\end{equation}
From the condition (\ref{eq:14}) and (9) follows condition (\ref{eq:13}).

{\bf Remark 2.} {\it It should be noted that theorem 2 is hold when the nonlinearities $Z_{1}, Z_{2}$ are continuous in the neighborhood of generating solution $z_{0}(n, c^{0})$.}

Now, we propose the following change of variables:
$$
z(n, \varepsilon) = z_{0}(n, c^{0}) + u(n, \varepsilon),
$$
where the element $c^{0}$ satisfies the operator equation (\ref{eq:14}). Then we can rewrite the boundary value problem (\ref{eq:11}), (\ref{eq:12}) in the following form
\begin{equation} \label{eq:15}
u(n + 1, \varepsilon) = A_{n}u(n, \varepsilon)  +\varepsilon \{ Z(z_{0}^{T}(n, c^{0}), n, 0) + Z^{'}_{u}(z_{0}^{T}(n, c^{0}), n, 0)u(n, \varepsilon) + \mathcal{R}(u(n, \varepsilon), n, \varepsilon)\},
\end{equation}
\begin{equation} \label{eq:16}
lu(\cdot, \varepsilon) = 0.
\end{equation}
Here $Z^{'}_{u}$ is Frechet derivative,
$$
\mathcal{R}(0, 0, 0) = \mathcal{R}^{'}_{u}(0, 0, 0) = 0.
$$
Boundary value problem (\ref{eq:15}), (\ref{eq:16}) has solutions if and only if the following condition is true:
\begin{equation} \label{eq:17}
\mathcal{P}_{N(\overline{Q}^{*})}l\sum_{i = 0}^{\cdot} \Phi(\cdot, i)(Z(z_{0}^{T}(i, c^{0}), i, 0) + Z_{u}^{'}(z_{0}^{T}(i, c^{0}), i, 0)u(i, \varepsilon) + \mathcal{R}(u(i, \varepsilon), i, \varepsilon)) = 0.
\end{equation}
Under this condition the set of solutions of boundary value problem (\ref{eq:15}), (\ref{eq:16}) has the following form
\begin{equation} \label{eq:18}
u(n, \varepsilon) = \mathcal{P}_{N(\overline{Q})}c + \overline{u}(n, \varepsilon),
\end{equation}
where
\begin{equation} \label{eq:19}
\overline{u}(n, \varepsilon) = \varepsilon \overline{G[Z(z_{0}^{T}(\cdot, c^{0}), \cdot, 0) + Z_{u}^{'}(z_{0}^{T}(\cdot, c^{0}), \cdot, 0)u(\cdot, \varepsilon) + \mathcal{R}(u(\cdot, \varepsilon), \cdot, \varepsilon), 0]}(n).
\end{equation}
Substituting (\ref{eq:18}) in (\ref{eq:17}) we obtain the following operator equation
\begin{equation} \label{eq:20}
B_{0}c = r,
\end{equation}
where the operator
$$
B_{0} = -\mathcal{P}_{N(\overline{Q}^{*})}l\sum_{i = 0}^{\cdot}\Phi(\cdot, i)Z_{u}^{'}(z_{0}^{T}(i, c^{0}), i, 0)\mathcal{P}_{N(\overline{Q})},
$$
$$
r = \mathcal{P}_{N(\overline{Q}^{*})}l\sum_{i = 0}^{\cdot}\Phi(\cdot, i)(Z_{u}^{'}(z_{0}^{T}(i, c^{0}), i, 0)\overline{u}(i, \varepsilon) + \mathcal{R}(u(i, \varepsilon), i, \varepsilon)).
$$
Using a generalization of implicit function theorem \cite{BoiPok} we have the following assertion.

{\bf Theorem 3.} (sufficient condition). {\it Suppose that the following condition is true:

 $\mathcal{P}_{N(\overline{B}_{0}^{*})}\mathcal{P}_{N(\overline{Q}^{*})} = 0.$

Then the boundary value problem (\ref{eq:11}), (\ref{eq:12}) has generalized solutions which can be found with using of iterative processes:
$$
u_{k + 1}(n, \varepsilon) = \mathcal{P}_{N(\overline{Q})}c_{k} + \overline{u}_{k}(n, \varepsilon),
$$
$$
c_{k + 1} = \overline{B}_{0}^{+}\mathcal{P}_{N(\overline{Q}^{*})}l\sum_{i = 0}^{\cdot}\Phi(\cdot, i)(Z_{u}^{'}(z_{0}^{T}(i, c^{0}), i, 0)\overline{u}_{k}(i, \varepsilon) + \mathcal{R}(u_{k}(i, \varepsilon), i ,\varepsilon)),
$$
$$
\overline{u}_{k + 1}(n, \varepsilon) = \varepsilon \overline{G[Z(z_{0}^{T}(\cdot, c^{0}), \cdot, 0) + Z_{u}^{'}(z_{0}^{T}(\cdot, c^{0}), \cdot, 0)\overline{u}_{k}(\cdot, \varepsilon) + \mathcal{R}(u_{k}(\cdot, \varepsilon), \cdot, \varepsilon), 0]}(n),
$$
where
$$
\mathcal{R}(u_{k}(n, \varepsilon), n, \varepsilon) = Z(z_{0}^{T}(n, c^{0}) + u_{k}(n, \varepsilon), n, \varepsilon) -
$$
$$
- Z(z_{0}^{T}(n, c^{0}), n, 0) - Z_{u}^{'}(z_{0}^{T}(n, c^{0}), n, 0)u_{k}(n, \varepsilon),
$$
$$
u_{0} = c_{0} = \overline{y}_{0} = 0.
$$
}

{\bf Applications.} It is well-known that systems like a Lotka-Volterra \cite{Volter}, \cite{Volt} plays an important role in the theoretical population biology \cite{Murray}, \cite{Murray1}. It is very important in mathematical biology as a model which describes dynamics of populations.  There exist many papers which are dedicated to investigation of such problems in continuous and discrete cases (see for example recently works \cite{RCL} - \cite{YWHWSR}). As a rule such problems are regular. We consider some examples of systems with different type of boundary conditions in the critical case. We show that the operator which generates considering problem can be Fredholm.  We find bifurcation conditions of solutions with using of the equation for generating constants.

{\bf Example 1.} Consider the following periodic boundary value problem in the finite dimensional case:
\begin{equation} \label{eq:21}
x_{i}(n + 1, \varepsilon) = a_{i}(n)x_{i}(n, \varepsilon) + b_{i}(n)y_{i}(n, \varepsilon) + \varepsilon g_{i}^{1}(n)x_{i}(n, \varepsilon)(1 - \sum_{j = 1}^{t}a_{ij}(n)y_{j}(n, \varepsilon)) + f_{1}^{i}(n),
\end{equation}
\begin{equation} \label{eq:22}
y_{i}(n + 1, \varepsilon) = c_{i}(n)x_{i}(n, \varepsilon) + d_{i}(n)y_{i}(n, \varepsilon) + \varepsilon g_{i}^{2}(n)y_{i}(n, \varepsilon)(1 - \sum_{j = 1}^{t}b_{ij}(n)x_{j}(n, \varepsilon)) + f_{2}^{i}(n),
\end{equation}
\begin{equation} \label{eq:23}
x_{i}(0, \varepsilon) = x_{i}(m, \varepsilon),
\end{equation}
\begin{equation} \label{eq:24}
y_{i}(0, \varepsilon) = y_{i}(m, \varepsilon),  i = \overline{1, p}.
\end{equation}
Here $x_{i}(n, \varepsilon), y_{i}(n, \varepsilon)$, $a_{i}(n), b_{i}(n), c_{i}(n), d_{i}(n), g_{i}^{1}(n), g_{i}^{2}(n), a_{ij}(n), b_{ij}(n) \in \mathbb{R}, i = \overline{1, p}, j = \overline{1, t}$.

For $\varepsilon = 0$ we obtain the following generating boundary value problem
\begin{equation} \label{eq:25}
x_{i}^{0}(n + 1) = a_{i}(n)x_{i}^{0}(n) + b_{i}(n)y_{i}^{0}(n) + f_{1}^{i}(n),
\end{equation}
\begin{equation} \label{eq:26}
y_{i}^{0}(n + 1) = c_{i}(n)x_{i}^{0}(n) + d_{i}(n)y_{i}^{0}(n) + f_{2}^{i}(n),
\end{equation}
\begin{equation} \label{eq:27}
x_{i}^{0}(0) = x_{i}^{0}(m),
\end{equation}
\begin{equation} \label{eq:28}
y_{i}^{0}(0) = y_{i}^{0}(m).
\end{equation}
In this case we have that the operator $l$ of boundary conditions has the following form:
$$l\left( \begin{array}{ccc} x_{0}(\cdot) \\
 y_{0}(\cdot) \end{array}\right) = \left( \begin{array}{ccc} x_{i}^{0}(m) - x_{i}^{0}(0) \\
  y_{i}^{0}(m) - y_{i}^{0}(0) \end{array} \right)_{i = \overline{1, k}} = \left( \begin{array}{ccc} 0 \\
  0 \end{array}\right).
  $$
For the vector $z_{i}^{0}(n) = (x_{i}^{0}(n), y_{i}^{0}(n))^{T}$ we can write the following assertion.

{\bf Corollary 1.} {\it The boundary value problem (\ref{eq:25})-(\ref{eq:28}) has periodic solutions if and only if
\begin{equation} \label{eq:29}
\mathcal{P}_{N(Q^{*})_{d}}\sum_{k = 0}^{m} \Phi(m, k)f(k) = 0,
\end{equation}
where $Q = \Phi(m, 0) - I$, $d$ is a number of linearly independent columns of $Q$;

under condition (\ref{eq:29}) the set of solutions has the form
\begin{equation} \label{eq:30}
z_{i}^{0}(n, c_{r}) = (G[f, 0])(n) + \mathcal{P}_{N(Q)_{r}}c_{r}, ~~ c_{r} \in \mathbb{R}^{r},
\end{equation}
where the generalized Green's operator $(G[f, 0])(n)$ has the following form
$$
(G[f, 0])(n) = -\Phi(n, 0)Q^{+}\sum_{k = 0}^{m}\Phi(m, k)f(k),
$$
$r$ is a number of linearly independent rows of $Q$.
}

{\bf Remark 3.} {\it It should be noted that in the considering case index of an operator $\mathcal{S}$ can be calculated in the following way
$$
ind~\mathcal{S} = r-d,
$$
where the operator $\mathcal{S}$ with boundary conditions has the following form
$$
\mathcal{S} \left( \begin{array}{ccc} x_{i}^{0}(n) \\
y_{i}^{0}(n) \end{array} \right): = \left( \begin{array}{ccc}  x_{i}^{0}(n + 1)  - a_{i}(n)x_{i}^{0}(n) - b_{i}(n)y_{i}^{0}(n)\\
y_{i}^{0}(n + 1) - c_{i}(n)x_{i}^{0}(n) - d_{i}(n)y_{i}^{0}(n)
\end{array}
\right).
$$
It means that the operator $\mathcal{S}$ is Fredholm \cite{BoiSam}.
}

For the nonlinear boundary value problem (\ref{eq:21})-(\ref{eq:24}) we obtain the following assertions.

{\bf Corollary 2.} (necessary condition). {\it If the boundary value problem (\ref{eq:21})-(\ref{eq:24}) has solution, then the element
$c_{r} = c_{r}^{0}$ satisfies the following equation for generating constants:
$$
F(c_{r}) = \mathcal{P}_{N(Q^{*})_{d}}\sum_{i = 0}^{m}\Phi(m, i)Z(z_{0}^{T}(i, c_{r}), i, 0) = 0,
$$
where
$$
Z(z_{0}^{T}(n, c_{r}), n, 0) = \left( \begin{array}{ccc} g_{i}^{1}(n)x_{i}^{0}(n, c_{r})(1 - \sum_{j = 1}^{t}a_{ij}(n)y_{j}^{0}(n, c_{r})) \\
g_{i}^{2}(n)y_{i}^{0}(n, c_{r})(1 - \sum_{j = 1}^{t}b_{ij}(n)x_{j}^{0}(n, c_{r})).
\end{array} \right)
$$
}

{\bf Corollary 3.} (sufficient condition). {\it Suppose that the following condition is true:

1) $\mathcal{P}_{N(B_{0}^{*})}\mathcal{P}_{N(Q^{*})_{d}} = 0$.

Then  the boundary value problem (\ref{eq:21})-(\ref{eq:24}) has generalized solution which can be found with using of iterative processes:
$$
u_{k + 1}(n, \varepsilon) = \mathcal{P}_{N(Q)_{r}}c_{k} + \overline{u}_{k}(n, \varepsilon),
$$
$$
c_{k + 1} = B_{0}^{+}\mathcal{P}_{N(Q^{*})}\sum_{i = 0}^{m}\Phi(m, i)(Z_{u}^{'}(z_{0}^{T}(i, c^{0}), i, 0)\overline{u}_{k}(i, \varepsilon) + \mathcal{R}(u_{k}(i, \varepsilon), i, \varepsilon)),
$$
$$
\overline{u}_{k + 1}(n, \varepsilon) = \varepsilon G[Z(z_{0}^{T}(\cdot, c^{0}), \cdot, 0) + Z_{u}^{'}(z_{0}^{T}(\cdot, c^{0}), \cdot, 0)\overline{u}_{k}(\cdot, \varepsilon) + \mathcal{R}(u_{k}(\cdot, \varepsilon), \cdot, \varepsilon), 0](n),
$$
where
$$
\mathcal{R}(u_{k}(n, \varepsilon), n, \varepsilon) = Z(z_{0}^{T}(n, c^{0}) + u_{k}(n, \varepsilon), n, \varepsilon) -
$$
$$
- Z(z_{0}^{T}(n, c^{0}), n, 0) - Z_{u}^{'}(z_{0}^{T}(n, c^{0}), n, 0)u_{k}(n, \varepsilon),
$$
$$
u_{0} = c_{0} = \overline{y}_{0} = 0,
$$
$$
Z_{u}^{'}(z_{0}^{T}(n, c^{0}), n, 0)u_{k}(n, \varepsilon) =
$$
$$
= \left(\begin{array}{ccc}
g_{i}^{1}(n)x_{i}^{0}(n, c_{r}^{0})(1 - \sum_{j = 1}^{t}a_{ij}(n)u_{jk}^{2}(n)) + g_{i}^{1}(n)u_{ik}^{1}(n)(1 - \sum_{j = 1}^{t}a_{ij}(n)y_{j}^{0}(n, c_{r}^{0})) \\
g_{i}^{2}(n)y_{i}^{0}(n, c_{r}^{0})(1 - \sum_{j = 1}^{t}b_{ij}(n)u_{jk}^{1}(n)) + g_{i}^{2}(n)u_{ik}^{2}(n)(1 - \sum_{j = 1}^{t}b_{ij}(n)x_{j}^{0}(n, c_{r}^{0}))
\end{array} \right).
$$
}
Suppose that $a_{i}(n) = b_{i}(n) = c_{i}(n) = g_{i}^{1}(n) = g_{i}^{2}(n) = a_{ij}(n) = b_{ij}(n) = 1$, $d_{i}(n) = 0$. In this case
$$A_{n} = A = \left( \begin{array}{ccc} 1 & 1 \\
1 & 0 \end{array} \right), n \in \mathbb{N}.$$ Then for the linear
boundary value problem (\ref{eq:25})-(\ref{eq:28}) we obtain that the evolution operator $\Phi(m, n)$ has the following form
$$
\Phi(m, n) = A^{m - n + 1} = \left( \begin{array}{ccc} F_{m - n + 2} & F_{m - n + 1} \\
 F_{m - n + 1} & F_{m - n} \end{array} \right).
$$
Here $F_{0} = 1, F_{1} = 1$, $F_{n + 2} = F_{n} + F_{n + 1}, n \geq 0$ are Fibonacci numbers.
In this case the matrix $Q$ is nondegenrate ($Q^{+} = Q^{-1}$, $\mathcal{P}_{N(Q)} = I, \mathcal{P}_{N(Q^{*})} = I$, $I$ is an identity matrix) and we obtain the following corollary.

{\bf Corollary 4.} {\it The boundary value problem (\ref{eq:21})-(\ref{eq:24}) has periodic solution if and only if
\begin{equation} \label{eq:100}
\sum_{k = 0}^{m}A^{m - k + 1}f(k) = \sum_{k = 0}^{m}\left( \begin{array}{ccc} F_{m - k + 2} & F_{m - k + 1} \\
F_{m - k + 1} & F_{m - k}   \end{array} \right) \left( \begin{array}{ccc} f_{1}^{i}(k) \\
f_{2}^{i}(k)\end{array} \right) = 0;
\end{equation}
under condition (\ref{eq:100}) the solution of the boundary value problem (\ref{eq:21})-(\ref{eq:24}) has the form
$$
z_{i}^{0}(n) = (G[f, 0])(n) = - A^{n + 1}Q^{-1}\sum_{k = 0}^{m}A^{m - k + 1}f(k)=
$$
$$
= -\frac{1}{\Delta(m)} \sum_{k = 0}^{m} \left( \begin{array}{ccc} a_{11}(n, m, k)f_{1}^{i}(k) + a_{12}(n, m, k)f_{2}^{i}(k) \\
a_{21}(n, m, k)f_{1}^{i}(k) +  a_{22}(n, m, k)f_{2}^{i}(k) \end{array} \right),
$$
where
$$
\Delta(m) = (F_{m + 2} - 1)(F_{m} - 1) - F_{m + 1}^{2};
$$
$$
a_{11}(n, m, k) = F_{n + 2} (F_{m}F_{m - k + 2}  - F_{m +1}F_{m - k +1}) - (F_{n + 2}F_{m - k + 2} + F_{n + 1}F_{m - k + 1}) +
$$
$$
+ F_{n + 1}(F_{m + 2}F_{m - k + 1} - F_{m + 1}F_{m - k + 2});
$$
$$
a_{12}(n, m, k) = F_{n + 2} (F_{m}F_{m - k + 1}  - F_{m +1}F_{m - k}) - (F_{n + 2}F_{m - k + 1} + F_{n + 1}F_{m - k}) +
$$
$$
+ F_{n + 1}(F_{m + 2}F_{m - k} - F_{m + 1}F_{m - k + 1});
$$
$$
a_{21}(n, m, k) = F_{n + 1} (F_{m}F_{m - k + 2}  - F_{m +1}F_{m - k +1}) - (F_{n + 1}F_{m - k + 2} + F_{n + 1}F_{m - k + 1}) +
$$
$$
+ F_{n}(F_{m + 2}F_{m - k + 1} - F_{m + 1}F_{m - k + 2});
$$
$$
a_{22}(n, m, k) = F_{n + 1} (F_{m}F_{m - k + 1}  - F_{m +1}F_{m - k}) - (F_{n + 2}F_{m - k + 1} + F_{n + 1}F_{m - k}) +
$$
$$
+ F_{n + 1}(F_{m + 2}F_{m - k} - F_{m + 1}F_{m - k + 1}).
$$
}
In this case the necessary condition of solvability for the nonlinear boundary value problem (\ref{eq:21})-(\ref{eq:24}) has the form
$$
\sum_{i = 0}^{m} \left( \begin{array}{ccc} F_{m - i + 2}x_{i}^{0}(n)(1 - \sum_{j = 1}^{n}y_{j}^{0}(n)) + F_{m - i + 1}y_{i}^{0}(n)(1 - \sum_{j = 1}^{n}x_{j}^{0}(n)) \\
F_{m - i + 1}x_{i}^{0}(n)(1 - \sum_{j = 1}^{n}y_{j}^{0}(n)) + F_{m - i}y_{i}^{0}(n)(1 - \sum_{j = 1}^{n}x_{j}^{0}(n)) \end{array} \right) = 0.
$$
In sufficient condition we have that
$$
Z_{u}^{'}(z_{0}^{T}(n), n, 0)u_{k}(n, \varepsilon) = \left(\begin{array}{ccc}
x_{i}^{0}(n)(1 - \sum_{j = 1}^{n}u_{jk}^{2}(n, \varepsilon)) + u_{ik}^{1}(n, \varepsilon)(1 - \sum_{j = 1}^{n}y_{j}^{0}(n)) \\
y_{i}^{0}(n)(1 - \sum_{j = 1}^{n}u_{jk}^{1}(n, \varepsilon)) + u_{ik}^{2}(n, \varepsilon)(1 - \sum_{j = 1}^{n}x_{j}^{0}(n))
\end{array} \right).
$$

{\bf Example 2.} Consider the following boundary value problem
\begin{equation} \label{eq:30}
x_{i}(n + 1, \varepsilon) = a_{i}(n)x_{i}(n, \varepsilon) + b_{i}(n)y_{i}(n, \varepsilon) + \varepsilon g_{i}^{1}(n)x_{i}(n, \varepsilon)(1 - \sum_{j = 1}^{n}a_{ij}(n)y_{j}(n, \varepsilon)) + f_{1}^{i}(n),
\end{equation}
\begin{equation} \label{eq:31}
y_{i}(n + 1, \varepsilon) = c_{i}(n)x_{i}(n, \varepsilon) + d_{i}(n)y_{i}(n, \varepsilon) + \varepsilon g_{i}^{2}(n)y_{i}(n, \varepsilon)(1 - \sum_{j = 1}^{n}b_{ij}(n)x_{j}(n, \varepsilon)) + f_{2}^{i}(n),
\end{equation}
with the following boundary conditions
\begin{equation} \label{eq:32}
l \left( \begin{array}{ccc} x_{i}(\cdot, \varepsilon) \\
y_{i}(\cdot, \varepsilon)\end{array}\right) = \left(\begin{array}{ccc} \sum_{k = 0}^{p_{1}}x_{i}(n_{k}, \varepsilon) \\ \sum_{l = 0}^{p_{2}}y_{i}(n_{l}, \varepsilon) \end{array}\right)_{i = \overline{1, p}} = \left(\begin{array}{ccc} \alpha_{1} \\
\alpha_{2} \end{array} \right).
\end{equation}
Here $n_{k}, k = \overline{0, p_{1}}$, $n_{l}, l = \overline{0, p_{2}}$ are finite sequences of integer numbers.
In this case we obtain the multi-point boundary-value problem.
Suppose that $x_{i}(n), y_{i}(n) \geq 0$ and boundary condition has the following form
\begin{equation} \label{eq:32}
l \left( \begin{array}{ccc} x_{i}(\cdot, \varepsilon) \\
y_{i}(\cdot, \varepsilon)\end{array}\right) = \left(\begin{array}{ccc} \sum_{i = 0}^{p}x_{i}(0, \varepsilon) \\ \sum_{i = 0}^{p}y_{i}(0, \varepsilon) \end{array}\right) = \left(\begin{array}{ccc} 1 \\
1 \end{array} \right).
\end{equation}
Such condition has practical meaning. It means the population distribution at the initial time (the proporion of the population in species).

\end{document}